\documentclass[11pt, oneside]{article}

\usepackage{graphicx}
\usepackage{amssymb}
\usepackage{amsmath}
\usepackage{amsthm}

\newcommand\norm[1]{\left\lVert#1\right\rVert}
\newcommand\floor[1]{\lfloor#1\rfloor}

\def\IN{{\mathbb N}}

\def\IR{{\mathbb R}}
\def\P{{\rm Prob}}
\def\E{{\bf E}}
\def\X{\mathcal{X}}
\def\L{\mathcal{L}}
\def\Heads{{\rm Heads}}
\def\Tails{{\rm Tails}}
\def\Leb{{\rm Leb}}
\def\half{{1 \over 2}}

\begin{document}

\centerline{\Large\bf The Coupling/Minorization/Drift Approach to}
\medskip
\centerline{\Large\bf Markov Chain Convergence Rates}

\bigskip \centerline{by (in alphabetical order)}

\medskip \centerline{Yu Hang Jiang, Tong Liu, Zhiya Lou, Jeffrey S. Rosenthal,
Shanshan Shangguan, Fei Wang, and Zixuan Wu}

\medskip \centerline{\sl Department of Statistical Sciences, University of Toronto}

\medskip\centerline{(August 18, 2020; last revised \today)}

\begin{quote}
\bigskip\noindent\textbf{Abstract}:
This review paper provides an introduction of Markov chains and their
convergence rates -- an important and interesting mathematical topic which
also has important applications for very widely used Markov chain Monte
Carlo (MCMC) algorithm.  We first discuss eigenvalue analysis for Markov
chains on finite state spaces.  Then, using the coupling construction, we
prove two quantitative bounds based on minorization condition and drift
conditions, and provide descriptive and intuitive examples to showcase how
these theorems can be implemented in practice. This paper is meant to
provide a general overview of the subject and spark interest in new Markov
chain research areas.
\end{quote}
\bigskip

\section{Introduction}

\begin{figure}[h]
\includegraphics[scale=0.1]{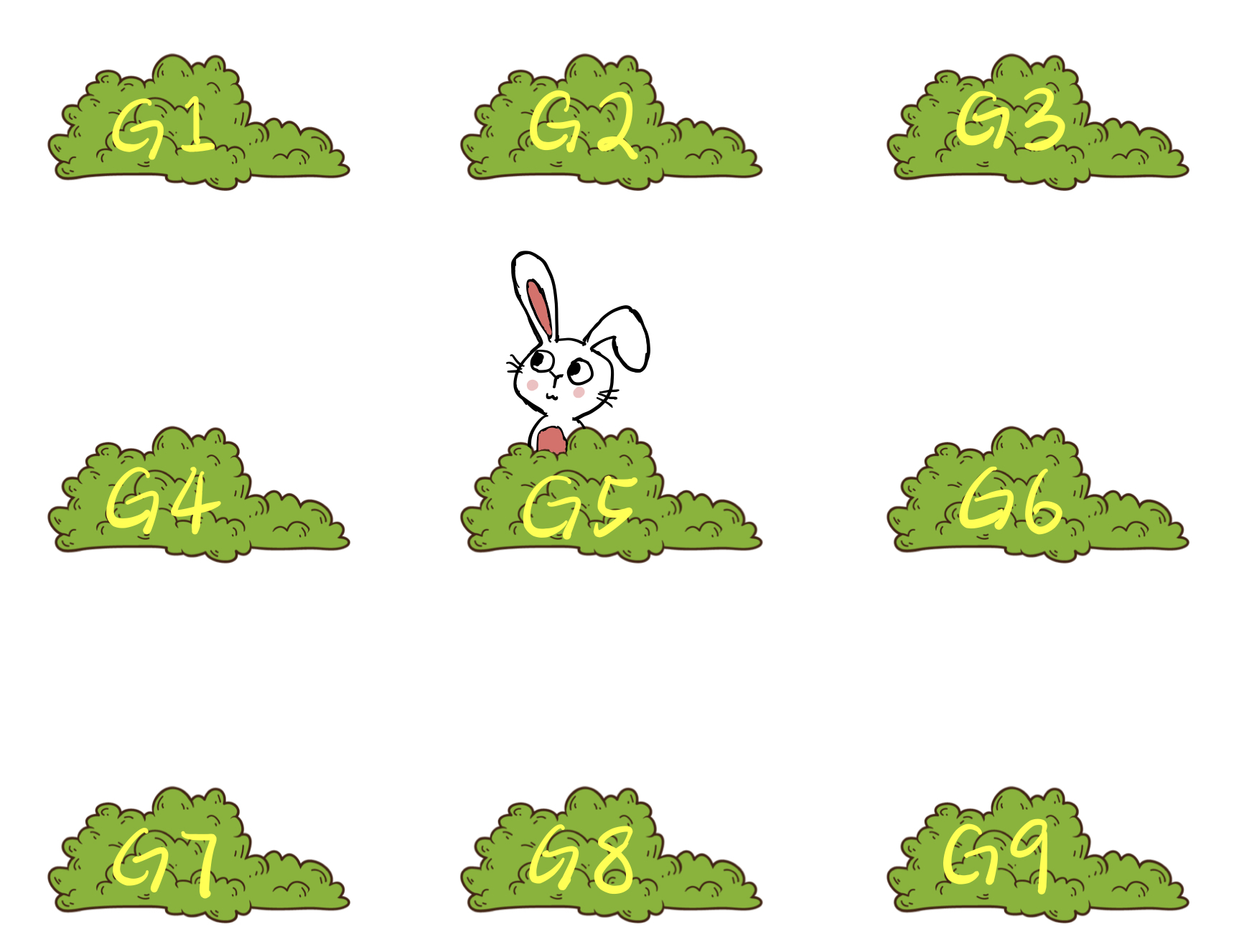}
\centering
\caption{The bunny example}
\end{figure}

\noindent
Imagine there is a $3\times 3$ grid of bushes, labeled $G_1, G_2,...,
G_9$, from top to bottom and left to right (Figure~1).
There is a fluffy little bunny
hiding in the middle bush, starving and ready to munch on some grass
around it. Assume the bunny never gets full and the grass is never
depleted. Once each minute, the bunny jumps from its current bush
to one of the nearest other bushes
(up, down, left, or right, not diagonal)
or stays at its current location, each with equal probability.
We can then ask about longer-term probabilities.
For example, if the bunny starts at G5, the probability of jumping to $G_7$
after two steps is:
\begin{align*}
    \P_{G5}(\mathrm{at}\  G_7 \ \mathrm{after\  two \ steps}) &=
\P(G_5 \to G_4 \to G_7) + \P(G_5 \to G_8 \to G_7)\\ &=
\frac{1}{5} \times \frac{1}{4} + \frac{1}{5} \times \frac{1}{4} =
\frac{1}{10}
\, .
\end{align*}
But what happens to the probabilities after three steps? ten steps? more?

This paper investigates the convergence of such
probabilities as the number of steps gets larger.
As we will discuss later, such bounds are not only an interesting topic in
their own right, they are also very important for reliably using Markov
chain Monte Carlo (MCMC) computer algorithms~\cite{metropolis,geyer,handbook}
which are very widely applied to numerous problems in statistics,
finance, computer science, physics, combinatorics, and more.
After reviewing the standard eigenvalue approach in Section~3, we will
concentrate on the use of ``coupling'', and specifically on the use of
``minorization'' (Section~4) and ``drift'' (Section~6) conditions.
We note that coupling is a very broad topic
with many different variations and applications (see
e.g.~\cite{thorisson}), and has even inspired its own algorithms (such as
``coupling from the past'').  And, there are many other methods of
bounding convergence of Markov chains, including continuous-time limits,
different metrics, path coupling, non-Markovian couplings,
spectral analysis, operator theory, and
more, as well as numerous other related topics,
which we are not able to cover here.

\section{Markov Chains}

The above bunny model is an example of a {\it Markov chain}
(in discrete time and space).
In general, a Markov chain is specified by three ingredients:

\textbf{1.} A {\it state space} $\mathcal{X}$, which is a collection of all of
the states the Markov chain might be at. In the bunny example, $\mathcal{X} =
\{G_1, G_2,..., G_9\}$.

\textbf{2.} An {\it initial distribution} (probability measure)
$\mu_0(\cdot)$, where $\mu_0(A)$ is the probability of starting within $A
\subset \mathcal{X}$ at time~0.  In the bunny example, $\mu_0(G_5)=1$, and
$\mu_0(G_i)=0\ \forall i\neq 5$.

\textbf{3.} A collection of {\it transition probability
distributions} $P(x, \cdot)$ on $\mathcal{X}$ for each state $x \in
\mathcal{X}$. The distribution $P(x,\cdot)$
represents the probabilities of the Markov
chain going from $x$ to the next state after one unit of time. In a
discrete state space like the bunny example, the transition probabilities
can be simply written as $P = \{p_{ij}\}_{i,j\in \mathcal{X}}$,
where $p_{ij}$ is the probability of jumping to $j$ from $i$.
Indeed, in the bunny example:

\centerline{\includegraphics[scale=0.15]{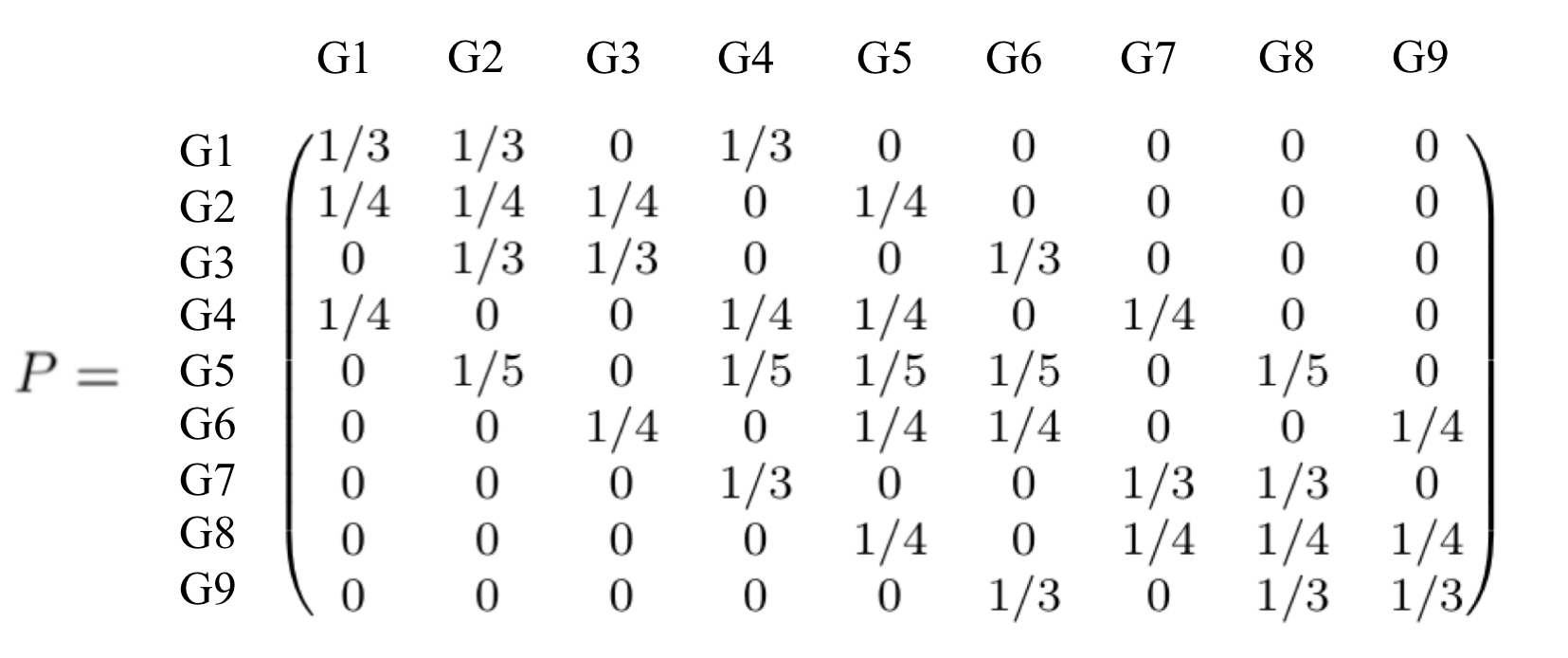}}

\noindent For example, in the second row, \(p_{21} = p_{22} = p_{23} =
p_{25} = \frac{1}{4}\) because from $G_2$, the probabilities of jumping to
each of $G_1$, $G_2$, $G_3$, or $G_5$ are each $\frac{1}{4}$.

We write $\mu_n(i)$ for the probability that the Markov chain is at
state~$i$ after~$n$ steps.  Given the initial distribution $\mu_0$ and
transition probabilities $P(x, \cdot)$, we can compute $\mu_n$ inductively by
$$
\mu_{n}(A) := \int_{x \in \mathcal{X}} P(x, A) \ \mu_{n-1}(dx)
\, ,
\qquad n \ge 1
\, .
$$
On a discrete space, this formula reduces to
$\mu_n(j) = \sum_{i \in \mathcal{X}} p_{ij} \, \mu_{n-1}(i)$.
In matrix form, regarding the $\mu_n$ as row-vectors, this means $\mu_n
= \mu_{n-1} \, P$. It follows by induction that $\mu_n = \mu_0  \, P^n$,
where $P^n$ is the $n$'th matrix power of $P$,
also called the $n$-step transition matrix.
Here $(P^n)_{ij}$ is the probability of jumping to $j$ from $i$
in $n$ steps. Indeed, if we take $\mu_0 = (0, 0, \dots, 1, \dots, 0)$,
so $\mu_0(i)=1$ with $\mu_0(j)=0$ for all $j\not=i$,
then $\mu_n(j) = \sum_{r=0} \mu_0(r) (P^n)_{rj}  = (P^n)_{ij}$.
This makes sense since if we start at $i$, then $\mu_n(j)$ is the
probability of moving from $i$ to $j$ in $n$ steps.

One main question in Markov chain analysis is whether the probabilities
$\mu_n$ will {\it converge} to a certain distribution, i.e.\ whether
$\pi := \lim_{n\to\infty} \mu_n$ exists.
If it does, then letting $n\to\infty$ in the relation
$\mu_{n+1} = \mu_n P$ indicates that $\pi$ must be {\it stationary},
i.e.\ $\pi = \pi P$.
On a finite state space, this means that $\pi$ is a left {\it eigenvector}
of the matrix $P$ with corresponding eigenvalue~1.

In the bunny example, by solving the system of linear equations
given by $\pi P = \pi$, the stationary probability distribution can be
computed to be the following vector:
$$
\pi \ = \
\bigg( \frac{1}{11}, \frac{4}{33}, \frac{1}{11}, \frac{4}{33}, \frac{5}{33}, \frac{4}{33}, \frac{1}{11}, \frac{4}{33}, \frac{1}{11} \bigg)
\, .
$$
In fact, the bunny example satisfies general theoretical properties called
\textit{irreducibility} and \textit{aperiodicity}, which guarantee that
the stationary distribution $\pi$ is unique, and that
$\mu_n$ converges to $\pi$ as $n\to\infty$ (see e.g.~\cite{spbook}).
However, in this paper we shall focus on {\it quantitative}
convergence rates, i.e.\ how large $n$ has to be to make $\mu_n$
sufficiently close to $\pi$.

\section{Eigenvalue Analysis on Finite State Spaces}
\label{sec-eigen}

When the state space is finite and small, it is sometimes possible to
obtain a quantitative bounds on the convergence rate through direct matrix
analysis (e.g.~\cite{eigen}).  We require eigenvalues $\lambda_i$ and
left-eigenvectors $v_i$ such that $v_i P = \lambda_i v_i$.  For example,
for the above bunny process, we compute (numerically, for simplicity) that
the eigenvalues and left-eigenvectors are:

\scriptsize
 $$(\lambda_0, \lambda_1, \lambda_2, \lambda_3, \lambda_4, \lambda_5, \lambda_6, \lambda_7, \lambda_8 ) = (1, 0.702, 0.702, -0.467, 0.333, 0.25, 0.25, 0.119, 0.119)$$
\[
\begin{pmatrix}
v_0\\ v_1\\ v_2\\ v_3\\ v_4\\ v_5\\ v_6\\ v_7\\ v_8\\ \end{pmatrix}
 = \begin{pmatrix}
0.091 & 0.121 & 0.091 & 0.121 & 0.152 & 0.121 & 0.091 & 0.121 & 0.091\\
0.489 & 0.361 & 0 & 0.361 & 0 & -0.361 & 0 & -0.361 & -0.489\\
0.055 & -0.318 & -4.863 & 0.399 & 0 & -0.399 & 4.863 & 0.318 & -0.055\\
-0.224 & 0.358 & -0.224 & 0.358 & -0.537 & 0.358 & -0.224 & 0.358 & -0.224\\
 0.500 & 0 & -0.500 & 0 & 0 & 0 & -0.500 & 0 & 0.500\\
 0.002 & -0.500 & 0.002 & 0.499 & -0.007 & 0.499 & 0.002 & -0.500 & 0.002\\
0.256 & -0.043 & 0.256 & -0.043 & -0.854 & -0.043 &  0.256 & -0.043 & 0.256\\
-0.436 & 0.394 & 0 & 0.394 & 0 & -0.394 & 0 & -0.394 & 0.436\\
0.018 & 0.377 & -0.425 & -0.410 & 0 & 0.410 & 0.435 & -0.377 & -0.018\\
\end{pmatrix}\\ \]\\ \normalsize

To be specific, assume that the bunny starts from the center
bush $G_5$, so $\mu_0 = (0,0,0,0,1,0,0,0,0)$.
We can express this $\mu_0$ in terms of the above eigenvector basis
as the linear combination:
$$
\mu_0 \ = \ v_0 -0.4255v_3 -0.7259v_6
\, .
$$
(Here $v_0 = \pi$, corresponding to the eigenvalue $\lambda_0=1$.)
Recalling that $\mu_n = \mu_0 P^n$, and that $v_i P = \lambda_i \, v_i$
by definition, we compute that e.g.\
$$
\mu_n(G_5)=(\lambda_0)^n v_0(G_5)-0.4255(\lambda_3)^n
v_3(G_5)-0.7259(\lambda_6)^n v_6(G_5)
$$
$$
= \pi(G_5)-0.4255(-0.4667)^n (-0.537)-0.7259(0.25)^n (-0.854)
\, .
$$
Since $|0.25| < |0.4667|$,
and $|0.4255 \cdot  (-0.537)| + |0.7259 \cdot (-0.854)| < 0.85$,
the triangle inequality implies that
$$
|\mu_n(G_5)-\pi(G_5)|
< 0.85 \, (0.4667)^n
\, ,
\quad n\in\IN
\, .
$$
This shows that $\mu_n(G_5) \to \pi(G_5)$, and gives a strong bound on the
difference between them.  For example,
$|\mu_n(G_5)-\pi(G_5)| < 0.01$ whenever $n \ge 6$, i.e.\
only 6 steps are required to make the bunny's probability of being at $G_5$
within 0.01 of its limiting (stationary) probability.  Other states
besides $G_5$ can be handled similarly.

Unfortunately, such direct eigenvalue or spectral analysis becomes
more and more challenging on larger and more complicated examples,
especially on non-finite state spaces.  So, we next introduce a different
technique which, while less tight, is more widely applicable.

\section{Coupling and Minorization Conditions}
\label{sec-minor}

The idea of {\it coupling} is to create two different copies of a random
object, and compare them.
Coupling has a long history in probability theory,
with many different applications and approaches (see e.g.~\cite{thorisson}).
A key idea is the {\it coupling inequality}.
Suppose we have two random variables $X$ and $Y$, each with their own
distribution.  Then for any subset $A$, we can write
\begin{multline*}
\big|\P(X\in A) - \P(Y \in A)\big| \\
\ = \
\big|\P(X\in A, \ X=Y) + \P(X\in A, \ X\not=Y) \\
- \P(Y\in A, \ X=Y) - \P(Y\in A, \ X\not=Y)\big|
\, .
\end{multline*}
But here $\P(X\in A, \ X=Y) = \P(Y\in A, \ X=Y)$, since they both
refer to the same event, so those two terms cancel.  Also, each of
$\P(X\in A, \ X\not=Y)$ and
$\P(Y\in A, \ X\not=Y)$ are between 0 and $\P(X\not=Y)$, so their
difference must be $\le \P(X\not=Y)$.  Hence,
$$
\big|\P(X\in A) - \P(Y \in A)\big| \ \le \ \P(X\not=Y)
\, .
$$
Since this upper bound is uniform over subsets $A$, we can even take a
supremum over~$A$, to also bound the {\it total variation distance}:
$$
\norm{\mathcal{L}(X) - \mathcal{L}(Y)}_{TV}
\ := \ \sup_{A\subseteq \mathcal{X}} \big|\P(X \in A) - \P(Y \in A)\big|
\ \le \ \P(X\not=Y)
\, .
$$
That is, the total variation distance between the probability laws
$\L(X)$ and $\L(Y)$ is bounded above by the probability that the random
variables $X$ and $Y$ are not equal.
To apply this fact to Markov chains,
the following condition is very helpful.

\paragraph{Definition.} A Markov chain with state space
$\mathcal{X}$ and transition probabilities
$P$ satisfies a \textit{minorization
condition} if there exists a (measurable)
subset $C \subseteq \mathcal{X}$, a probability
measure $\nu$ on $\mathcal{X}$, a constant $\epsilon > 0$, and a positive
integer $n_0$, such that
$$
P^{n_0}(x, \cdot) \geq \epsilon \nu(\cdot), \ x \in C
\, .
$$
We call such $C$ a \textit{small set}. In particular, if $C =
\mathcal{X}$ (the entire state space), then the Markov chain
satisfies a \textit{uniform minorization condition}, also referred to
as {\it Doeblin's condition}.
\bigskip

For a concrete example, suppose the state space is the half-line
$\mathcal{X}=[0,\infty)$, with transition probabilities given by
\begin{equation}\label{concreteex}
P(x,dy) \ = \
\Big( e^{-2y}+\frac{1}{\sqrt{2\pi}(x+1)}e^{-\frac{y^2}{2(x+1)^2}} \Big) \, dy
\, .
\end{equation}
That is, from a state $x$, the chain moves to an equal mixture of an
Exponential(2) distribution and a half-normal distribution with
mean~0 and standard deviation $x+1$.
In this case, $P(x,dy) \ge e^{-2y} \, dy$ for all $x$ (see Figure~2),
so the chain satisfies a
uniform minorization condition with $n_0=1$, $\nu(y) = 2e^{-2y}$,
and $\epsilon = \frac{1}{2}$.

\begin{figure}[h]
\includegraphics[scale=0.15]{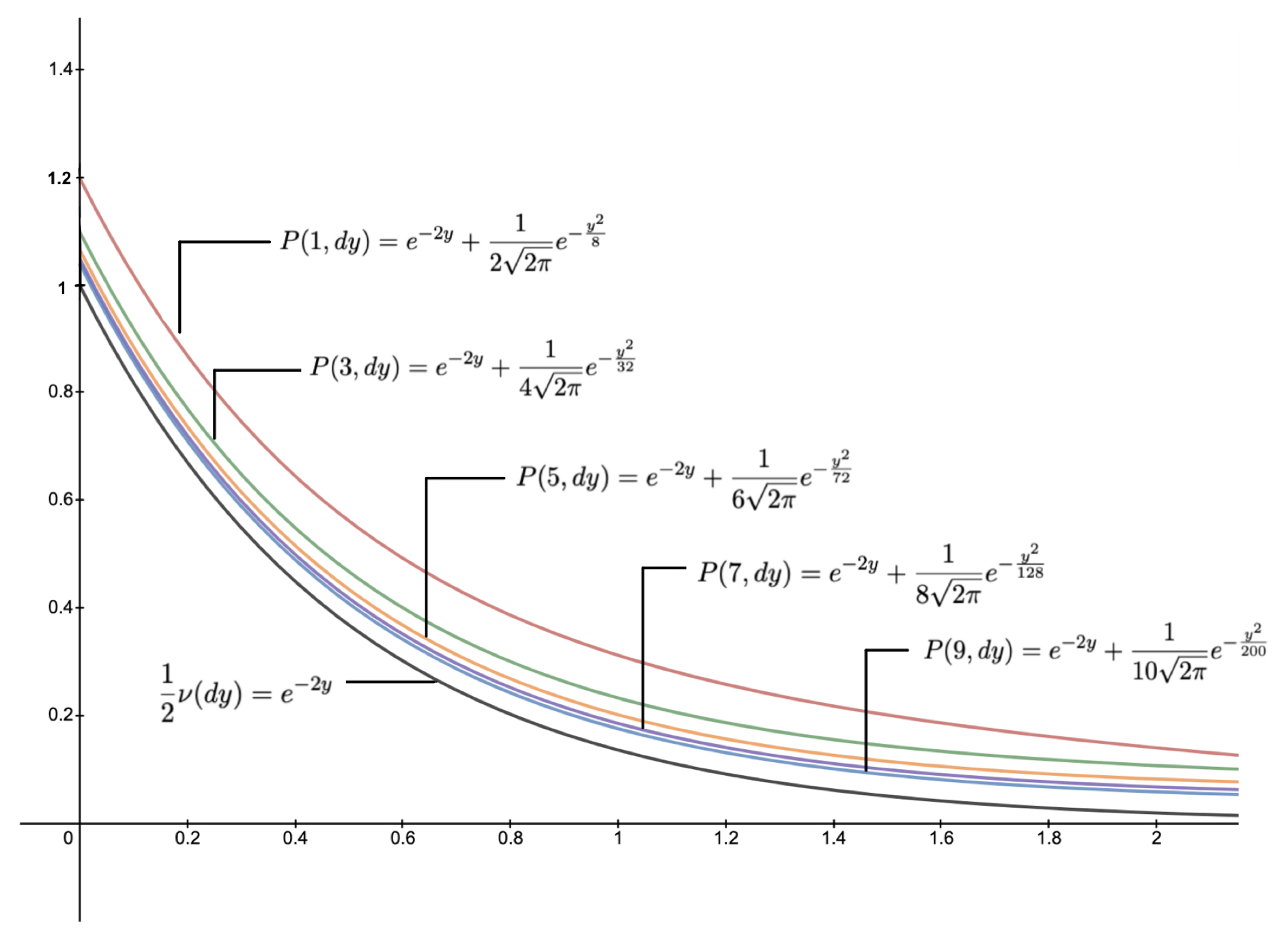}
\centering
\caption{The example satisfying the minorization condition.}
\end{figure}

The uniform minorization condition implies that there exists a
common overlap of size $\epsilon$
between all of the transition probabilities.
This allows us to formulate a coupling construction of two
different copies $\{X_n\}$ and $\{X'_n\}$ of a Markov chain, as
follows.  Assume for now that $n_0=1$.
First, choose $X_0 \sim \mu_0(\cdot)$ and $X_0' \sim \pi(\cdot)$ independently.
Then, inductively for $n=0,1,2,\ldots$,

\textbf{1.} If $X_n = X_n'$, choose $z \sim P(X_n, \cdot)$ and let
$X'_{n + 1} = X_{n + 1} = z$. The chains have already coupled, and they will
remain equal forever.

\textbf{2.} If $X_n \not= X_n'$, flip a coin whose probability of Heads
is $\epsilon$.  If it shows Heads,
choose $z \sim \nu(\cdot)$ and let $X'_{n + {1}} = X_{n + {1}} = z$.
Otherwise, update $X_{n + {1}}$ and $X'_{n + 1}$ independently with
probabilities given by
$$
\P(X_{n + 1} \in A) = \frac{P(X_n, A) -
\epsilon \nu(A)}{1 - \epsilon},
\quad \P(X'_{n + 1} \in A) = \frac{P(X'_n, A) -
\epsilon \nu(A)}{1 - \epsilon}
\, .
$$
(The minorization condition guarantees that these ``residual'' probabilities
are non-negative, and hence are probability measures since their total
mass equals~${P(X_n,\X) - \epsilon \nu(\X) \over 1-\epsilon}
= {1-\epsilon \over 1-\epsilon} = 1$.)
This construction ensures that overall, $\P(X_{n+1}\in A |X_n = x) = P(x,A)$
and $\P(X'_{n+1}\in A |X'_n = x)= P(x,A)$ for any $x \in
\mathcal{X}$: indeed, if the two chains are unequal at time $n$, then
\begin{equation*}
\begin{split}
& \P(X_{n + 1} \in A \ | \ X_n = x) \\
&= \P(X_{n + 1} \in A, \, \Heads \ | \
X_n = x ) +\P(X_{n + 1} \in A,  \, \Tails \ | \ X_n = x) \\
&=\P(\Heads) \ \P(X_{n + 1} \in A \ | \ X_n = x, \ \Heads) \\
&\qquad\qquad\qquad + \P(\Tails)
\ \P(X_{n + 1} \in A \ | \ X_n = x, \ \Tails) \\
&= \epsilon \, \nu(A) + (1 - \epsilon) \, \frac{P(x, A) -
\epsilon \, \nu(A)}{1 - \epsilon}
\ = \ P(x, A)
\, .
\end{split}
\end{equation*}

If $n_0 > 1$, then we can use the above construction for the times
$n=0,n_0,2n_0,\ldots$, with $n+1$ replaced by $n+n_0$,
and with $P(\cdot,\cdot)$ replaced by
$P^{n_0}(\cdot,\cdot)$.  Then, if desired, we can later
``fill in'' the intermediate states $X_n$ for $jn_0 < n < (j+1)n_0$,
from their appropriate conditional distributions
given the already-constructed values of $X_{jn_0}$ and $X_{(j+1)n_0}$.

Now, since $X'_0 \sim \pi(\cdot)$, and $\pi$ is a stationary distribution,
therefore $X'_n \sim \pi(\cdot)$ for all $n$.
And, every \(n_0\) steps, the two chains
probability at least
\(\epsilon\) of coupling (i.e., of the coin showing Heads).
So, $\P(X_n \not= X_{n}') \le (1 - \epsilon)^{\lfloor  n/ n_0 \rfloor }$,
where $\lfloor \cdot \rfloor$ means {\it floor}.
The coupling equality then implies:

\paragraph{Theorem~1.} If $\{X_n\}$ is
a Markov chain on $\mathcal{X}$, whose transition probabilities
satisfy a uniform minorization
condition for some $\epsilon > 0$,
then for any positive integer $n$, and any $x\in\X$,
$$
\norm{\mathcal{L}(X_n) - \pi(\cdot)}_{TV}
\ \le \ (1 - \epsilon)^{\lfloor  n/ n_0 \rfloor }
.
$$
\smallskip

For the above Markov chain~\eqref{concreteex}, we showed a uniform
minorization condition with $n_0=1$ and $\epsilon=1/2$.  So, Theorem~1
immediately implies that
$\norm{\mathcal{L}(X_n) - \pi(\cdot)}_{TV}
\le (1 - \epsilon)^{\lfloor  n/ n_0 \rfloor }
= (1 - (1/2))^n = 2^{-n}$,
which is $< 0.01$ if $n \ge 6$, i.e.\
this chain converges within 6 steps.

If $\X$ is finite, and for some $n_0\in\IN$
there is at least one state $j\in\X$ such that
the $j^{\rm th}$ column of $P^{n_0}$ is all positive, i.e.\
$(P^{n_0})_{ij} > 0$ for all $i\in\X$.
Then we can set
\(\epsilon = \sum_{j \in \mathcal{X}}\min_{i \in \mathcal{X}} (P^{n_0})_{ij}
> 0\), and
$\nu(j) = \epsilon^{-1} \, \min_{i \in \mathcal{X}} (P^{n_0})_{ij}$,
so that $(P^{n_0})_{ij} \geq \epsilon \, \nu(j)$ for all $i,j \in \mathcal{X}$,
i.e.\ an \(n_0\)-step uniform minorization condition is satisfied with that
value of $\epsilon$.

\subsection{Application to Bunny Example}

The bunny example does not satisfy a one-step minorization condition,
since every column of $P$ has some zeroes,
so we instead consider its two-step transition probabilities, $P^2$:

\centerline{ \includegraphics[scale=0.15]{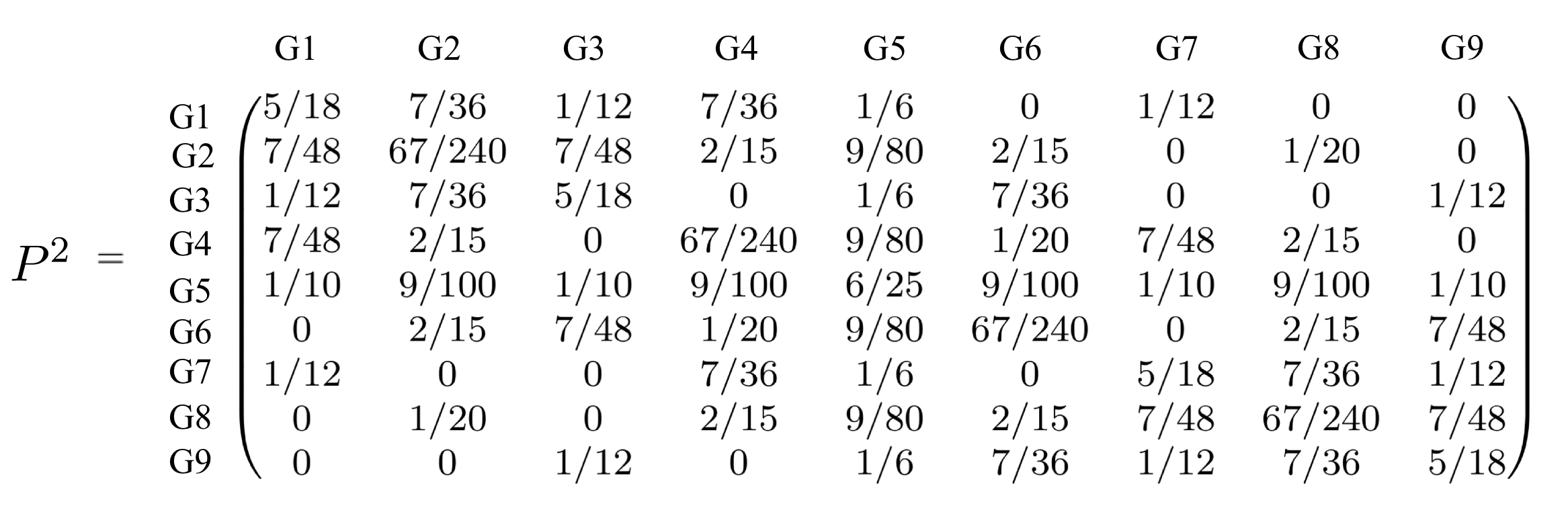} }

\noindent
In this two-step transition matrix,
the fifth column only contains positive values, since
no matter where the bunny starts, there
will always be at least a $\frac{9}{80}$
chance that it will jump to the center bush (\(G5\)) in two steps.
Thus, we can satisfy a two-step minorization condition by taking
$$
    \epsilon = \sum_{j \in \mathcal{X}} \min_{i \in \mathcal{X}} (P^2)_{ij}
            = 0 + \dots + 0 + \frac{9}{80} + 0 + \dots + 0 \\
            = \frac{9}{80}
$$
and $\nu(j) = \epsilon^{-1} \, \min_{i \in \mathcal{X}} (P^{n_0})_{ij}$
as above.  Then, we can apply Theorem~1, with $n_0=2$ and
$\epsilon = 9/80$, to conclude that
\[
\norm{\mathcal{L}(X_n) - \pi(\cdot)}_{TV}
\ \le \ \Big(1-\frac{9}{80}\Big)^{\floor{n/2}}
= \ \Big(\frac{71}{80}\Big)^{\floor{n/2}}
\]
For example, if we want the
distribution of the bunny's location to be within 0.01 of the
stationary distribution $\pi$, this is achieved within $n=78$ steps.
This bound is not nearly as tight as our previous
result $n=6$, but it is uniform over all states (not just $G5$),
plus it was derived using a much more general method (without
the need to compute eigenvalues and eigenvectors).
Of course, such bounds might be more difficult to obtain
on larger, more complicated examples.

\subsection{Pseudo-Minorization Conditions}

The coupling construction used to prove Theorem~1 was a {\it pairwise}
construction, i.e.\ it only considered two chain locations $x$
and $y$ at a time.  If we replace the distribution \(\nu(\cdot)\)
with \(\nu_{xy}(\cdot)\), allowing it to depend on \(x\) and \(y\),
then \(C\) is called a \textit{pseudo-small set}, and
Theorem~1 continues to hold~\cite{pseudosmall}.
It then follows that on a finite state space, if we instead choose
$$
\epsilon
\ = \ \min_{i,j \in \mathcal{X}} \sum_{z \in \mathcal{X}}
\min\{(P^{n_0})_{iz}, (P^{n_0})_{jz}\} > 0;
\quad \nu_{ij}(z) = \frac{\min \{(P^{n_0})_{iz}, (P^{n_0})_{jz}\}}
{\sum\limits_{w \in \mathcal{X}}\min \{(P^{n_0})_{iw}, (P^{n_0})_{jw}\}}
\, ,
$$
then the chain will satisfy an \(n_0\)-step pseudo-minorization condition,
i.e.\ for all $i,j,z \in \mathcal{X}$,
$(P^{n_0})_{iz} \geq \epsilon \, \nu_{ij}(z)$ and
$(P^{n_0})_{jz} \geq \epsilon \, \nu_{ij}(z)$.
Hence, exactly as above, we will again have
$\norm{\mathcal{L}(X_n) - \pi(\cdot)}_{TV}
\, \le \, (1 - \epsilon)^{\lfloor  n/ n_0 \rfloor}$.

We now apply this pseudo-minorization idea
to the bunny example, with $n_0=2$.  Examining the
matrix $P^2$ above, we see that the minimum values of \(\sum_{z \in
\mathcal{X}} \min \{(P^2)_{iz}, (P^2)_{jz}\}\) occur at
\((i,j) = (3, 7)\) or \((1, 9)\), corresponding to
opposite corners of the \(3 \times 3\) grid (which
makes sense since opposite corners will have the least amount of
transitional overlap).  We then calculate the minorization
constant
$$
\epsilon
\ = \ \sum_{z \in \mathcal{X}} \min \{(P^2)_{3z}, (P^2)_{7z}\}
\ = \ \frac{1}{12} + 0 + 0 + 0 + \frac{1}{6} + 0 + 0 +
0 + \frac{1}{12} \ = \ \frac{1}{3}
\, .
$$
Therefore,
$\norm{\mathcal{L}(X_n) - \pi(\cdot)}_{TV}
\, \le \, (1-\epsilon)^{\floor{n/2}}
\, = \, (\frac{2}{3})^{\floor{n/2}}$.
For instance, this bound is $< 0.01$ if $n=24$, i.e.\ if
the bunny jumps 24 times. This is a significant improvement
over the previous minorization result of $n=78$,
though it is still not as tight as the specific eigenvalue bound of $n=6$.

\section{Continuous State Space: Point Process MCMC}

The above analysis was primarily focused on finite state spaces, such as
the bunny example.  We now extend to continuous examples on subsets of
$\IR^d$.

To be specific, consider a {\it point process} consisting of
three particles each randomly located within the closed rectangle
$[0, 1]^2 \subset \mathbb{R}^2$, with positions
denoted by $x = (x_i)_{i = 1, 2, 3} =
(x_{i1}, x_{i2})_{i = 1, 2, 3}$, so the state space $\mathcal{X} = [0,
1]^6$. Suppose these particles are distributed according to
a probability distribution with {\it unnormalized density}
(meaning that the actual
density is a constant multiple of so it integrates to 1) given by
$$
\pi(x)
\ := \ \pi(x_1,x_2,x_3)
\ = \ \exp\Big[ - C\sum_{i = 1}^3||x_i||
			- D\sum_{i < j} ||x_i - x_j||^{-1} \Big]
\, ,
$$
where $C$ and $D$ are fixed positive constants, and $||\cdot||$ is the
usual Euclidean ($L^2$) norm on $\mathbb{R}^2$. In this density, the first
sum pushes the particles towards the origin, and the second sum pushes
them away from each other.

We now create a Markov chain which has $\pi$ as its stationary
distribution.  To do this, we use a version of the {\it Metropolis
Algorithm}~\cite{metropolis}.  Each step of the Markov chain proceeds
as follows.  Given $X_n=x$, we first ``propose'' to move the particles
from their current configuration $x$ to some other configuration $y$,
chosen from the uniform (i.e., Lebesgue) measure on $\X$.  Then, with
probability $\min[1, \, {\pi(y) \over \pi(x)}]$, we ``accept'' this
proposal and move to the new configuration by setting $X_{n+1} = y$.
Otherwise, we ``reject'' this proposal and leave the configuration
unchanged by setting $X_{n+1}=x$.

This Metropolis Algorithm is a well-known procedure which can easily be
shown~\cite{metropolis,handbook} to create a Markov chain which has
$\pi$ as its stationary distribution.  It is the most common type of
{\it Markov chain Monte Carlo} (MCMC) algorithm.  Such algorithms are a
very popular and general method of generating samples from complicated
probability distributions, by running the corresponding Markov chain
for many iterations.  They are used very frequently in a wide variety
of fields, ranging from Bayesian statistics to financial modeling to
medical research to machine learning and more.  For further background,
see e.g.~\cite{handbook} and the many references therein.

However, to get reliable samples, it is important to know how many
iterations are required to approximately converge to $\pi$, i.e.\ to
establish quantitative convergence bounds.  The above Markov chain has
an uncountably infinite state space $\X$, so the eigenvalue analysis
of Section~\ref{sec-eigen} is not easily available (though there have
been some efforts to use spectral analysis on general state spaces,
see e.g.\ \cite{geyer} and other papers).  On the other hand,
the uniform minorization condition of Section~\ref{sec-minor} can still
be applied.  Indeed, in the web appendix~\cite{appendix}, we prove:

\paragraph{Lemma~1.} The Markov chain constructed above satisfies a
uniform minorization condition with $n_0=1$ and
$\epsilon = (0.48) \, e^{-C (4.25) - D(9.88)}$.
\bigskip

For example, if $C = D = 1/10$, then we may take $\epsilon=0.117$.
It then follows from Theorem~1 that we have the convergence bound
$$
\norm{\mathcal{L}(X_n) - \pi(\cdot)}_{TV}
\ \le \ (1-\epsilon)^n
\ = \ (1- 0.117)^n
\ = \ (0.883)^n
\, .
$$
This shows that after $n=38$ steps, the total variation distance between
our Markov chain and the stationary distribution will be less than 0.01.

\section{Unbounded State Space: Drift Conditions}

In the previous section, the uniform minorization condition give us a good
quantitative convergence bound.  However, in many cases, especially on
unbounded state spaces, the minorization condition cannot be satisfied
uniformly, only on some subset $C \subseteq \X$.  In such cases, we have
to adjust our previous $n_0=1$ coupling construction, as follows.  We first
choose $X_0 \sim \mu(\cdot)$ and $X_0' \sim \pi(\cdot)$ independently,
and then inductively for $n=0,1,2,\ldots$,

\textbf{1}. If $X_n=X^{\prime}_n$, we choose $X_{n+1}=X^{\prime}_{n+1}
\sim P(X_n,\cdot)$.

\textbf{2}. Else, if $(X_n,X^{\prime}_n) \in C \times C$,
we flip a coin whose probability of Heads
is $\epsilon$, and then
update $X_{n + 1}$ and $X'_{n + 1}$ in the same way as in step~2
of our previous (uniform minorization) construction above.

\textbf{3}. Else,
if $(X_n,X^{\prime}_n) \not\in C \times C$, then we just
conditionally independently choose $X_{n+1} \sim
P(X_n,\cdot)$ and $X_{n+1}^{\prime} \sim P(X_n^{\prime},\cdot)$, i.e.\
the two chains are simply updated independently.

\medskip The above construction provides good coupling bounds provided
that the two chains return to $C \times C$ often enough, but this last
property is difficult to guarantee.  Thus, to obtain convergence bounds,
we also require a {\it drift condition}. Basically, the drift condition
guarantees that the chains will return to $C \times C$ quickly enough that
we can still achieve a coupling.

\bigskip\noindent
\textbf{Definition.} A Markov chain with a small set $C \subseteq \X$
satisfies a \textit{bivariate drift condition} if there exists a function
$h: \mathcal{X} \times \mathcal{X} \rightarrow [1, \infty)$ and some
$\alpha > 1$, such that
\begin{equation*}
    \Bar{P}h(x,y) \leq h(x,y) / \alpha, \quad (x,y) \notin C \times C
\, ,
\end{equation*}
where
$$
\Bar{P}h(x,y) \ := \
\E[h(X_{n+1},Y_{n+1}) \ | \ X_n = x, Y_n = y]
$$
is the expected (average) value of $h(X_{n+1},Y_{n+1})$ on the next
iteration, when the chains start from $x$ and $y$ respectively
(and proceed independently).

\bigskip Such
bivariate drift conditions can be combined with non-uniform
minorization conditions to produce quantitative convergence bounds.  To
state them, we use the quantity
$B_{n_0} = \max\{1,\alpha^{n_0} (1-\epsilon) \sup_{C \times C} \Bar{R}h\}$,
where
$$
\Bar{R}h(x, y) \ = \ \int_{\mathcal{X}} \int_{\mathcal{X}}
(1-\epsilon)^{-2} \, h(z,w) \,
[P^{n_0}(x,dz)-\epsilon \nu(dz)] \, [P^{n_0}(y,dw) - \epsilon \nu(dw)]
\, .
$$
This daunting expression represents the expected
value of \(h(X_{n+n_0}, X'_{n+n_0})\) given that \(X_n = x\), that
\(X'_n = y\), and that
the two chains fail to couple at time $n$ (i.e.\ the corresponding
coin shows Tails).  We can simplify $\Bar{R}h$ in certain situations.
For example, if \(D\) is a set such that \(P^{n_0}(x, D) = 1\) for all
\(x \in C\), then since the expected value of a random variable is
always less than the maximal value it could take, we have
$$
\sup_{(x,y) \in C \times C} \Bar{R}h
\ \le \ \sup_{(x, y) \in D \times D} h(x,y)
\, .
$$

With all that in mind, we have the following result:

\paragraph{Theorem~2.} Consider a Markov Chain on $\mathcal{X}$, with
$X_0=x$, and transition probabilities $P$. Suppose the above minorization
and bivariate drift conditions hold, for some $C \subseteq \mathcal{X}$,
$h: \mathcal{X} \times \mathcal{X} \rightarrow [1,\infty)$, probability
distribution $\nu(\cdot)$, $\alpha > 1$, and $\epsilon > 0$.
Then for any integers $1 \leq j \leq n$, with $B_{n_0}$ as above,
$$
\norm{\mathcal{L}(X_n) - \pi}_{TV}
\ \le \ (1-\epsilon)^j + \alpha^{-n} B_{n_0}^{j-1}
\E_{Z \sim \pi}[h(x,Z)]
\, .
$$
\smallskip

Here we give a basic idea of the proof; for more details,
see~\cite{computsimple,probsurv}.  We create a second copy of the Markov chain
with $X'_0 \sim \pi$, and use the above coupling construction.
Let $N_n$ be the number of times the chain
$(X_n,X_n^{\prime})$ is in $C \times C$ by the $n^{\rm th}$ step.
Then by the coupling inequality,
$$
\norm{\mathcal{L}(X_n) - \pi}_{TV}
\ \le \ P[X_n \neq X^{\prime}_n]
\ \le \ P[X_n \neq X_n^{\prime}, N_{n-1} \geq j] + P[X_n \neq X_n^{\prime},
N_{n-1} < j]
\, .
$$
The first term suggests
that the chains have not coupled by time $n$ despite visiting $C \times C$
at least $j$ times. Since each such time gives them a chance of $\epsilon$
to couple, the first term is $\le (1 - \epsilon)^j$. The
second term is more complicated, but from the bivariate drift
condition together with a martingale argument, it
can be shown to be no greater than $\alpha^{-n} B_{n_0}^{j-1}
\E_{Z\sim\pi}[h(x,Z)]$.

Sometimes it could be hard to directly check the bivariate drift
condition. We now introduce the more easily-verified univariate drift
condition, and give a way to derive the bivariate condition from the
univariate one.

\paragraph{Definition.} A Markov chain with a small set $C$
satisfies a \textit{univariate drift condition} if there are constants
$0< \lambda <1 $ and $b<\infty$, and a function $V:\
\mathcal{X}\rightarrow[1,\infty]$, such that:
$$
PV(x) \ \le \ \lambda V(x) +
b\mathbf{1}_C(x), \quad x \in \mathcal{X}
\, ,
$$
where $PV(x) := \E[V(X_{n + 1}) \ | \ X_n = x]$.

This univariate drift condition can be used to bound $\E_\pi(V)$.
Indeed, assuming $\E_\pi(V)<\infty$, stationarity then implies that
$\E_\pi(V) \le \E_\pi(V) + b$, whence $\E_\pi(V) \le b/(1-\lambda)$.
But also, univariate drift condition can imply bivariate drift conditions,
as follows.

\paragraph{Proposition.} Suppose the univariate drift condition is
satisfied for some $V:\ \mathcal{X}\rightarrow[1,\infty]$, $C \in
\mathcal{X}$, $0<\lambda <1 $ and $b<\infty$. Let $d= \inf_{x\in C^c}
V(x)$. If $d > [b/(1-\lambda)]-1$, then the bivariate drift condition is
satisfied for the same C, with $h(x,y)=\frac{1}{2}[V(x)+V(y)]$ and
$\alpha^{-1} = \lambda +b/(d+1) <1$.

\begin{proof}
Assume $(x,y) \notin C \times C$. Then either $x \notin C$ or $y \notin
C$, so $h(x,y) \ge (1+d)/2$. Then, our univariate drift condition
applied separately to $x$ and to $y$ implies that
$PV(x)+PV(y) \le \lambda V(x)+ \lambda V(y) +b$. Therefore
\begin{align*}
    \bar{P}h(x,y) &= \frac{1}{2}[PV(x)+PV(y)] \le \frac{1}{2}[\lambda
V(x)+ \lambda V(y) +b] = \lambda h(x,y) + b/2 \\ &\le \lambda h(x,y) +
(b/2)[h(x,y)/((1+d)/2)]
    =[\lambda + b/(1+d)] h(x,y)
\, ,
\end{align*}
which gives the result.
\end{proof}

We now apply these non-uniform quantitative convergence bounds to a Markov
chain on an unbounded state space.
Let the state space be $\X = \mathbb{R}$, the entire real line,
with unnormalized target density $\pi(x) = e^{-|x|}$.

To create a Markov chain which has $\pi$ as its stationary
distribution, we use another version of the Metropolis
Algorithm~\cite{metropolis,handbook}.
Each step of the Markov chain proceeds as follows.  First, we
propose to move from the state
$x$ to some other state $y$, chosen from the uniform (i.e.,
Lebesgue) measure on the interval $[x-2,x+2]$.
Then, with probability $\min[1, \ {\pi(y) \over \pi(x)}]$,
we accept this proposal and move to the new state $y$,
otherwise we reject it and remain at $x$.
Once again, this procedure
creates a Markov chain which has $\pi$ as its stationary distribution.

To apply Theorem~2, we need to establish minorization and drift conditions.
In the web appendix~\cite{appendix}, we prove:

\paragraph{Lemma~2.} The above Markov chain satisfies a minorization
condition with $C=[-2,2]$, $n_0=2$, $\epsilon = {1 \over 8 e^2}$, and
$\nu(A) = \half \, \Leb(A \cap [-1,1])$, where $\Leb$ is Lebesgue measure
on $\IR$.

\paragraph{Lemma~3.} The above Markov chain satisfies a univariate drift
condition with $V(x) = e^{-|x|/2}$,
$C = [-2, 2]$, $\lambda=0.916$, and $b=0.285$.
\bigskip

We can then apply the above Proposition to derive
a bivariate drift condition. Note that here
$d = \inf_{x\in C^c} V(x) = e$, and $[b/(1 - \lambda)] - 1 = 2.39 < e$.
So, $h(x,y) := \frac{1}{2}(V(x) + V(y))$ satisfies a bivariate drift
condition with $\alpha^{-1} = \lambda + b/(d+1)
= 0.916 + 0.285/(e + 1) \doteq 0.993$.

We also need to bound the above quantity $B_{n_0} = B_2$.
Let $D = [-6, 6]$.  Then clearly $P^2(x, D) = 1$ for any $x \in C$. Thus
$$
\sup_{(x, y) \in C \times C} \bar{R}h(x,y)
\ \le \ \sup_{(x, y) \in D \times D} h(x ,y)
\ = \ \sup_{x \in D}V(x) \ = \ e^3 \ < \ 20.1
\, .
$$
So $B_2 \equiv \max[1,\alpha^2 (1-\epsilon) \sup \bar{R}h]
< (0.993)^{-2} (1-{1 \over 8 e^2}) (20.1) \doteq 20.04$.

Let \(X_0 = 0\).  Then
\begin{multline*}
    \E_{Z\sim\pi}[h(0,Z)]
    \ = \ \E_{Z\sim\pi}[\frac{1}{2}(V(0) + V(Z))]
    \ = \ \frac{1}{2} + \half \E_\pi(V)
\\
    \ = \ \frac{1}{2} + \half \frac{
\int_{y \in \mathcal{X}} e^{\frac{1}{2}|y|}e^{-|y|}dy}
{\int_{y \in \mathcal{X}}e^{-|y|}dy}
    \ = \ \half + \half \times \frac{2}{1}
    \ = \ 2
\, .
\end{multline*}
(If this specific calculation were not available, then we could instead
use the bound $\E_\pi(V) \le b/(1-\lambda) = 0.285/(1-0.916) = 3.393$ as
discussed above.)
Therefore, by Theorem~2, with $X_0=0$, we have
\begin{multline*}
\norm{\mathcal{L}(X_n) - \pi}_{TV}
\ \le \ (1-\epsilon)^j + \alpha^{-n} B_2^{j-1} \E_{Z\sim\pi}[h(0, Z)]
\\
\ \leq \ (0.983)^j + (0.993)^n(20.04)^{j-1}[2]
\, .
\end{multline*}
For example, setting $n=120,000$ and $j=274 = 1 + n/439.56$, this becomes
$
\norm{\mathcal{L}(X_n) - \pi(\cdot))}_{TV}
\ \leq \ (0.983)^{274} + [(0.993)(20.04)^{1/439.56}]^{120000} [2]
\, < \, 0.01
\, ,
$
i.e.\ the Markov chain is within 0.01 of stationarity after 120,000
iterations.  This is quite a conservative upper bound.  Nevertheless, we
have obtained a concrete
quantitative convergence bound, for an unbounded Markov chain.

\section{Conclusion}

This paper has discussed Markov chains and their convergence rates, and
why they are important for MCMC algorithms.  We introduced the eigenvalue
method, the coupling method, and minorization and drift conditions, and
applied them to examples on state spaces ranging from finite to compact
to unbounded.  For the bunny example, we showed several possible methods
of obtaining convergence bounds. Indeed, bounding a Markov chain's
convergence rate is not an one-time, definitive process; for various
Markov chains, it is possible to strengthen the bound through careful
and creative new constructions. The bounds presented in this paper all
have their imperfections, and will certainly not give tight or realistic
bounds for all examples.  There is plenty of room for new and tighter
and more flexible bounds, which can help us to understand Markov chains
better, and also run MCMC algorithms more confidently and reliably.

\bigskip\noindent\bf Acknowledgements. \rm
We thank the editor and reviewers for very helpful comments on the first
version of this manuscript.

\end{document}